\documentclass[a4paper,12pt]{article}

\usepackage{t1enc}
\usepackage{amsmath}
\usepackage{amssymb}
\usepackage{amsfonts}

\usepackage[dvips]{graphicx}

\usepackage{epic}

\newtheorem{CONJ}{Conjecture}
\newtheorem{LEM}{Lemma}
\newenvironment{DEM}{\noindent {\bf Proof:}}{\hfill $\Box$}
\newcommand{\ds}{\displaystyle}

\title{Finding largest small polygons with GloptiPoly}

\begin{document}

%\begin{article}

%\begin{opening}

\author{Didier Henrion$^{1,2,3}$,
Fr\'ed\'eric Messine$^4$}

%\runningauthor{Didier Henrion, Fr\'ed\'eric Messine}
%\runningtitle{Finding largest small polygons with GloptiPoly}

\footnotetext[1]{CNRS; LAAS; 7 avenue du colonel Roche, F-31077 Toulouse; France. {\tt henrion@laas.fr}}
\footnotetext[2]{Universit\'e de Toulouse; UPS, INSA, INP, ISAE; UT1, UTM, LAAS; F-31077 Toulouse; France}
\footnotetext[3]{Faculty of Electrical Engineering, Czech Technical University in Prague,
Technick\'a 2, CZ-16626 Prague, Czech Republic}
\footnotetext[4]{ENSEEIHT-IRIT, 2 rue Charles Camichel, BP 7122, F-31071 Toulouse; France. {\tt messine@n7.fr}}

\maketitle

\begin{abstract}
A small polygon is a convex polygon of unit diameter. We are interested in small polygons which have the largest area for a given number of vertices $n$.
Many instances are already solved in the literature, namely for all odd $n$, and for $n=4, 6$ and $8$. Thus, for even $n\geq 10$, instances
of this problem remain open. Finding those largest small polygons can be formulated as nonconvex quadratic programming
problems which can challenge state-of-the-art global optimization algorithms. We show that a recently developed
technique for global polynomial optimization, based on a semidefinite programming
approach to the generalized problem of moments and implemented in the public-domain
Matlab package GloptiPoly, can successfully find largest small polygons for $n=10$ and $n=12$.
Therefore this significantly improves existing results in the domain. When coupled with
accurate convex conic solvers, GloptiPoly can provide numerical
guarantees of global optimality, as well as rigorous guarantees relying
on interval arithmetic.
\end{abstract}

\begin{center}
\end{center}

%\end{opening}

\section{Introduction}

The problem of finding the largest small polygons was first studied by Reinhardt in 1922 \cite{Rei1922}. He solved the problem by proving that the solution corresponds to the regular small polygons but only when the number of vertices $n$ is odd. He also solved the case $n=4$ by proving that a square with diagonal length equal to 1 is a solution. However, there exists an infinity of other solutions (it is just necessary that the two diagonals intersect with a right angle). The hexagonal case $n=6$ was solved numerically by Graham in 1975 \cite{Gra1975}. Indeed, he studied possible structures that the optimal solution must have. He introduced the diameter graph of a polygon which is defined by the same vertices as the polygon and by edges if and only if the corresponding two vertices of the edge are at distance one. Using a result due to Woodall \cite{Woo1971}, he proved that the diameter graph of the largest small polygons must be connected, yielding 10 distinct possible configurations for $n=6$. Discarding 9 of these 10 possibilities (by using standard geometrical reasonings plus the fact that all the candidates must have an area greater than the regular small hexagon), he determined the only possible diameter graph configuration which can provide a better solution than the regular one. He solved this last case numerically, yielding the largest small hexagon. The name of this corresponding optimal hexagon is \emph{Graham's little hexagon}. Following the same principle, Audet et al. in 2002 found the \emph{largest small octagon} \cite{Aud2002}. The case $n=8$ is much more complicated than the case $n=6$ because it generates 31 possible configurations and just a few of them can be easily discarded by geometrical reasonings. Furthermore, for the remaining cases, Audet et al. had to solve difficult global optimization problems with 10 variables and about 20 constraints. In \cite{Aud2002}, these problems are formulated as quadratic programs with quadratic constraints. For solving this program, Audet et al. used a global solver named QP \cite{Aud2000}. Notice that optimal solutions for $n=6$ and $n=8$ are not the regular polygons \cite{Aud2002, Gra1975}. At the Toulouse Global Optimization Workshop TOGO 2010, the corresponding optimal octagon was named \emph{Hansen's little octagon}. The following cases $n=10, 12, \ldots$ were open. However in 1975, Graham proposed a conjecture which is the following: when $n$ is even and $n\geq 4$, the largest small polygon must have a diameter graph with a cycle with $n-1$ vertices and with an additional edge attached to a vertex of the cycle; this is true for $n=4, 6$ and also $n=8$, see Figure~\ref{Fig-31} for an example of this optimal structure for the octagon. Therefore, this yields only one possible diameter graph configuration that must have the optimal shape. In 2007, Foster and Szabo proved Graham's conjecture \cite{Fos2007}. Thus to solve the following open cases $n \geq 10$ and $n$ even, it is just necessary to solve one global optimization problem defined by the configuration of the diameter graph with a cycle with $n-1$ vertices and an additional pending edge. In order to have an overview of these and other related subjects about polygons, refer to \cite{Aud2007,Aud2009}.

In this paper, using the Global Optimization software GloptiPoly \cite{gloptipoly3}, we solve the two open cases $n=10$ and $n=12$.  Invariance of the quadratic programming problem under a group of permutation suggests that small polygons have a symmetry axis. Exploiting this symmetry
allows to reduce signficantly the size of optimization problems to be solved, even though currently we are not able to prove that this is without loss of generality.
In Section~\ref{sec:QP}, the general quadratic formulation is presented. Then, we show that the quadratic problem is invariant under a group of permutations, which
suggests an important reduction of the size of the quadratic programs. In Section~\ref{sec:GloptiPoly} GloptiPoly is presented and an example of its use is given to solve the case of the octagon yielding a rigorous certificate
of the global optimum found in \cite{Aud2002}. In Section~\ref{sec:NUM}, we present the solutions for the largest small decagon ($n=10$) and dodecagon ($n=12$), and conjecture the solutions for the tetradecagon ($n=14$) and the hexadecagon ($n=16$).
Then, we conclude in Section~\ref{CONC}.

%\end{article}\end{document}

\section{Nonconvex Quadratic Optimization Problems}
\label{sec:QP}

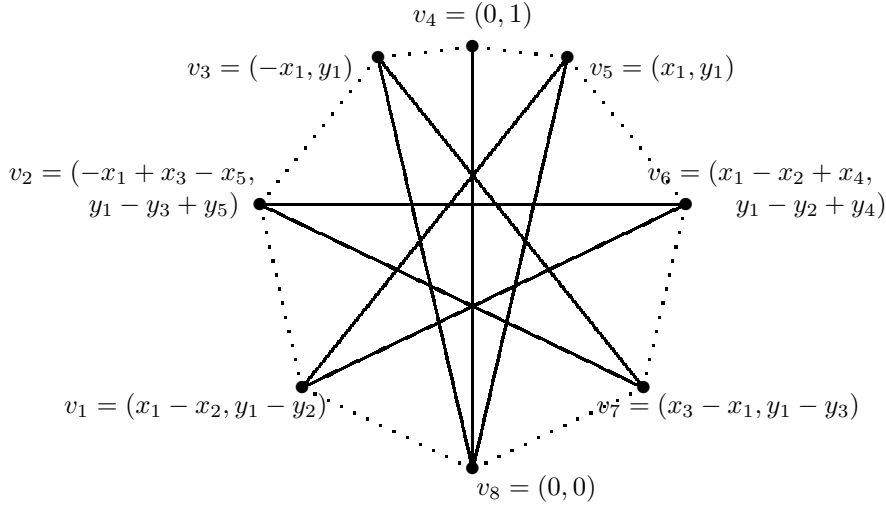
\begin{figure}[ht]
\setlength{\unitlength}{2.8mm}
\begin{picture}(40,22)(0,0)
\thicklines
 \dottedline{0.1}(20, 1)(20, 21)
 \dottedline{0.1}(20, 1)(24.450, 20.4985)
 \dottedline{0.1}(20, 1)(15.550, 20.4985)
 \dottedline{0.1}(15.550, 20.4985)(28.019, 4.862)
 \dottedline{0.1}(24.450, 20.4985)(11.981, 4.862)
 \dottedline{0.1}(10,     13.540)(28.019, 4.862)
 \dottedline{0.1}(30,     13.540)(11.981, 4.862)
 \dottedline{0.1}(30,     13.540)(10,     13.540)
 \put(20,1){\makebox(0,0){$\bullet$}}
 \put(20,21){\makebox(0,0){$\bullet$}}
 \put(24.45,20.50){\makebox(0,0){$\bullet$}}
 \put(15.55,20.50){\makebox(0,0){$\bullet$}}
 \put(28.02,4.86){\makebox(0,0){$\bullet$}}
 \put(11.98,4.86){\makebox(0,0){$\bullet$}}
 \put(30,13.54){\makebox(0,0){$\bullet$}}
 \put(10,13.54){\makebox(0,0){$\bullet$}}
\dottedline{0.7}(20,
1)(28.02,4.86)(30,13.54)(24.45,20.5)(20,21)(15.55,20.5)(10,13.54)(11.98,4.86)(20,1)
%  \put(-3 , 0){\vector(1, 0){42}}
%  \put(-3 , 0){\vector(0, 1){21}}
%  \put(40, 0){\makebox(0,0){$x$}}
%  \put(-3, 22){\makebox(0,0){$y$}}
  \put(23,  .1){\makebox(0,0){\footnotesize $v_8=(0,0)$}}
  \put( 7,   4){\makebox(0,0){\footnotesize $v_1=(x_1-x_2, y_1-y_2)$}}
  \put( 4,  15){\makebox(0,0){\footnotesize $v_2=(-x_1+x_3-x_5,$}}
  \put(5.5, 13.5){\makebox(0,0){\footnotesize $y_1-y_3+y_5)$}}
  \put(10.5,20){\makebox(0,0){\footnotesize $v_3=(-x_1, y_1)$}}
  \put(20,  22.5){\makebox(0,0){\footnotesize $v_4=(0,1)$}}
  \put(28.9,  20){\makebox(0,0){\footnotesize $v_5=(x_1, y_1)$}}
  \put(33.5,15){\makebox(0,0){\footnotesize $v_6=(x_1-x_2+x_4,$}}
  \put(35.9,13.5){\makebox(0,0){\footnotesize $y_1-y_2+y_4)$}}
  \put(32,   4){\makebox(0,0){\footnotesize $v_7=(x_3-x_1, y_1-y_3)$}}
\end{picture}
\caption{Case of $n=8$ vertices. Definition of variables following Graham's conjecture.} \label{Fig-31}
\end{figure}

As mentioned above, for even $n\geq 4$, finding the largest small polygon with $n$ vertices amounts to solving only one global optimization problem \cite{Fos2007}, namely a nonconvex quadratic programming problem. This formulation was previously introduced by Audet et al. in \cite{Aud2002} for the octagon case. On Figure~\ref{Fig-31}, we recall the nomenclature and configuration
corresponding to the octagon ($n=8$).

\begin{equation}
\left\{
\begin{array}{rll}
 \displaystyle \max_{x,y} & \multicolumn{2}{l}{A_8=x_1+\frac{1}{2}\left\{ (x_2+x_3-4x_1)y_1 + (3x_1-2x_3+x_5)y_2
                             \right.}\\
                   &  \multicolumn{2}{l}{\quad \left.+ (3x_1-2x_2+x_4)y_3 + (x_3-2x_1)y_4 + (x_2-2x_1)y_5 \right\}} \\
     \mbox{s.t.}   & \|v_i-v_j\|^2 \leq 1 & i,j=1,\ldots,8\\
                   & \|v_2-v_6\|^2 = 1\\
                   & x_i^2 + y_i^2 = 1     & i = 1,\ldots,5 \\
                   & y_i \geq 0  & i = 1,\ldots,5  \\
     & 0 \leq x_1 \leq \frac{1}{2} \\
  & 0 \leq x_i \leq 1 & i = 2,\ldots,5. \\
 \end{array}
 \right.
 \label{eq-31}
\end{equation}

Without loss of generality we can insert the additional constraint $x_2 \geq x_3$ which eliminates a symmetry axis. In program (\ref{eq-31}),
all the constraints are quadratic.
The quadratic objective function corresponds to the computation of the area of the octagon following Graham's diameter graph configuration.

\begin{figure}[ht]
\setlength{\unitlength}{0.9mm}
\begin{picture}(40,22)(0,0)
\thicklines
 \dottedline{0.1}(20, 1)(20, 21)
 \dottedline{0.1}(20, 1)(24.450, 20.4985)
 \dottedline{0.1}(20, 1)(15.550, 20.4985)
 \put(20,1){\makebox(0,0){$\bullet$}}
 \put(20,21){\makebox(0,0){$\bullet$}}
 \put(24.45,20.50){\makebox(0,0){$\bullet$}}
 \put(15.55,20.50){\makebox(0,0){$\bullet$}}
  \put(20, -1){\makebox(0,0){\footnotesize $(0,0)$}}
  \put(9,19){\makebox(0,0){\footnotesize $(-x_1, y_1)$}}
  \put(20,  23.5){\makebox(0,0){\footnotesize $(0,1)$}}
  \put(29.8,  19){\makebox(0,0){\footnotesize $(x_1, y_1)$}}
\end{picture}
\begin{picture}(40,22)(0,0)
\thicklines
 \dottedline{0.1}(24.450, 20.4985)(11.981, 4.862)
 \dottedline{0.1}(10,     13.540)(28.019, 4.862)
 \dottedline{0.1}(30,     13.540)(11.981, 4.862)
 \dottedline{0.1}(30,     13.540)(10,     13.540)
 \dottedline{0.7}(20, 1)(20, 21)
 \dottedline{0.7}(20, 1)(24.450, 20.4985)
 \put(20,1){\makebox(0,0){$\bullet$}}
 \put(20,21){\makebox(0,0){$\bullet$}}
 \put(24.45,20.50){\makebox(0,0){$\bullet$}}
 \put(28.02,4.86){\makebox(0,0){$\bullet$}}
 \put(11.98,4.86){\makebox(0,0){$\bullet$}}
 \put(30,13.54){\makebox(0,0){$\bullet$}}
 \put(10,13.54){\makebox(0,0){$\bullet$}}
  \put(28.9,  20){\makebox(0,0){\footnotesize $u_0$}}
  \put( 7,   4){\makebox(0,0){\footnotesize $u_1$}}
  \put(5.5, 13.5){\makebox(0,0){\footnotesize $u_3$}}
  \put(33.5,15){\makebox(0,0){\footnotesize $u_2$}}
  \put(32,   4){\makebox(0,0){\footnotesize $u_4$}}
  \put(20, -1){\makebox(0,0){\footnotesize $(0,0)$}}
  \put(20,  23.5){\makebox(0,0){\footnotesize $(0,1)$}}
\end{picture}
\begin{picture}(40,22)(0,0)
\thicklines
 \dottedline{0.7}(20, 1)(20, 21)
 \dottedline{0.7}(20, 1)(15.550, 20.4985)
 \dottedline{0.1}(15.550, 20.4985)(28.019, 4.862)
 \dottedline{0.1}(10,     13.540)(28.019, 4.862)
 \dottedline{0.1}(30,     13.540)(10,     13.540)
 \put(20,1){\makebox(0,0){$\bullet$}}
 \put(20,21){\makebox(0,0){$\bullet$}}
 \put(15.55,20.50){\makebox(0,0){$\bullet$}}
 \put(28.02,4.86){\makebox(0,0){$\bullet$}}
 \put(30,13.54){\makebox(0,0){$\bullet$}}
 \put(10,13.54){\makebox(0,0){$\bullet$}}
  \put( 4,  15){\makebox(0,0){\footnotesize $w_2$}}
  \put(10.5,20){\makebox(0,0){\footnotesize $w_0$}}
  \put(33.5,15){\makebox(0,0){\footnotesize $w_3$}}
  \put(32,   4){\makebox(0,0){\footnotesize $w_1$}}
  \put(20, -1){\makebox(0,0){\footnotesize $(0,0)$}}
  \put(20,  23.5){\makebox(0,0){\footnotesize $(0,1)$}}
\end{picture}
\medskip
\caption{Construction of the General Quadratic Formulation.} \label{Fig-3c}
\end{figure}
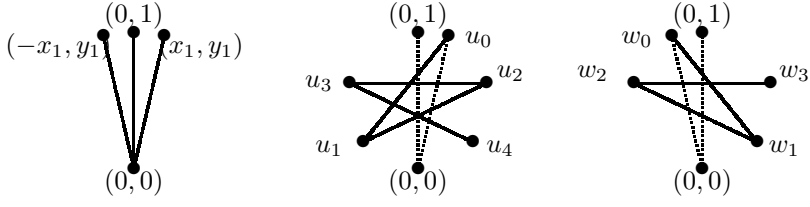

\medskip
To generalize this quadratic formulation, we have to define two vertices $v_n=(0,0)$ and $v_\frac n2=(0,1)$ and $n-2$ other vertices $v_i$ with the help
of the following variables:
\begin{eqnarray}
u_k&=&\left(x_1+\sum_{i=1}^k(-1)^ix_{2i}, y_1+\sum_{i=1}^k(-1)^iy_{2i}\right),\nonumber\\
w_k&=&\left(\sum_{i=0}^k(-1)^{i+1}x_{2i+1}, \sum_{i=0}^k(-1)^iy_{2i+1}\right),\nonumber
\end{eqnarray}
with $k=0,1,\ldots,\frac n2-2$. On Figure \ref{Fig-3c}, we represent the four first vertices, namely $(0,0), (-x_1,y_1), (0,1), (x_1,y_1)$. Note that we have an axis of symmetry from the line passing through $(0,0)$ and $(0,1)$, i.e. this symmetry provides the maximal area for this quadrilateral polygon and its corresponding value is equal to $x_1$. Thus, from vertex $u_0=(x_1,y_1)$, we construct iteratively all the vertices $u_i$ following a path in the graph of diameter and we do the same from $w_0=(-x_1,y_1)$ for the vertices $w_i$, see Figure \ref{Fig-3c}.

\bigskip

The main idea using this notation is that it yields:
$$\|u_{k+1}-u_k\|=\left\|\left((-1)^{k+1}x_{2(k+1)}, (-1)^{k+1}y_{2(k+1)}\right)\right\|=x_{2(k+1)}^2+y_{2(k+1)}^2=1.$$
We obtain by the same way that  $\|w_{k+1}-w_k\|=x_{2(k+1)+1}^2+y_{2(k+1)+1}^2=1$.
This allows to eliminate all the variables $y_i$, reducing by half the size of the problem.

\medskip
Then, we rename the vertices by $v_i=(\overline{x}_i,\overline{y}_i)$:
$$
\left\{
\begin{array}{ll}
v_i:=u_{2i-1}& i=1,\ldots,\lfloor\frac{n-2}{4}\rfloor\\
v_{n-i}:=w_{2i-1}&  i=1,\ldots,\lfloor\frac{n-2}{4}\rfloor\\
v_{\frac n2-i}:=w_{2(i-1)}& i=1,\cdots,\lceil\frac{n-2}{4}\rceil\\
v_{\frac n2+i}:=u_{2(i-1)}& i=1,\cdots,\lceil\frac{n-2}{4}\rceil\\
\end{array}
\right.
$$

The general quadratic program can be written as follows:
\begin{equation}
\left\{
\begin{array}{rll}
 \displaystyle \max_{x,y} & \multicolumn{2}{l}{A_n=x_1+\frac{1}{2}\ds\sum_{i=1\atop i\neq \frac n2-1, i\neq \frac n2}^{n-2}(\overline{y}_i\overline{x}_{i+1}-\overline{x}_i\overline{y}_{i+1})}\\
     \mbox{s.t.}   & \|v_i-v_j\|^2 \leq 1 & i,j=1,\ldots,n \\
                   & \|v_{\lfloor\frac n4\rfloor}-v_{\lceil\frac {3n}4\rceil}\|^2 = 1\\
                   & x_i^2 + y_i^2 = 1     & i=1,\ldots,n-3 \\
                   & y_i \geq 0 & i=1,\ldots,n-3 \\
                   & 0 \leq x_1 \leq \frac{1}{2}  \\
& 0 \leq x_i \leq 1 &  i=1,\ldots,n-3 \\
 \end{array}
 \right.
 \label{quad-prog}
\end{equation}
where $\overline{x}_i$ and $\overline{y}_i$ are linear functions respectively depending on $x_i$ and $y_i$.

We also give equivalent formulations for the computation of the area of an $n$-gon:
\begin{eqnarray*}
A_n&=&\frac{1}{2}\ds\sum_{i=1}^{n-2}(\overline{y}_i\overline{x}_{i+1}-\overline{x}_i\overline{y}_{i+1}),\\
A_n&=&x_1+\frac{1}{2}\left(\ds\sum_{i=1}^{\frac n2-2}(\overline{y}_i\overline{x}_{i+1}-\overline{x}_i\overline{y}_{i+1})+\ds\sum_{i=\frac n2+1}^{n-2}(\overline{y}_i\overline{x}_{i+1}-\overline{x}_i\overline{y}_{i+1})\right),\\
A_n&=&\frac{1}{2}\sum_{i=1}^n (\overline{x}_i+\overline{x}_{i+1})\times (\overline{y}_i-\overline{y}_{i+1}),
\end{eqnarray*}
where $i+1$ is taken modulus $n$ in the second expression of $A_n$.

\medskip
Remark that the largest little quadrilateral, Graham's little hexagon and Hansen's little octagon all have an axis of symmetry passing through vertices $(0,0)$ and $(0,1)$ in our quadratic formulation. In the following, we state a result suggesting that this symmetry axis should occur for all the largest small $n$-gons when $n$ is even.

Let us define the substitution $\sigma$ as follows:
$$
\sigma(x_{2i})=x_{2i+1}\mbox{ and }\sigma(x_{2i+1})=x_{2i}, \quad i=1,\ldots,n-3.
$$
This provides a group $G$ of substitutions which permutes all the vertices except $v_{\frac n2-1}$ and $v_{\frac n2+1}$ through the symmetry axis passing by $v_n$ and $v_\frac n2$.
\begin{LEM}\label{group}
Quadratic problem~(\ref{quad-prog}) is invariant through substitution group $G$.
\end{LEM}

\begin{DEM}
This substitution does not affect the constraints because it is just a new numbering of the vertices and $\sigma\left(v_{\lfloor\frac n4\rfloor}\right)=v_{\lceil\frac {3n}4\rceil}$ and reciprocally $\sigma\left(v_{\lceil\frac {3n}4\rceil}\right)=v_{\lfloor\frac n4\rfloor}$. Hence, it remains to prove that $A_n$ is invariant by $\sigma$.

Denote by $\hat v_i=(\hat x_i, \hat y_i)$ with $i=1,\ldots n$ the new vertices obtained using substitution $\sigma$. Hence, $\hat v_i=(-\overline{x}_{n-i},\overline{y}_{n-i})$, for all $i=1,\ldots, \frac n2-2$ and the other vertices are unchanged, i.e, $\hat v_n=v_n, \hat v_\frac n2=v_\frac n2, \hat v_{\frac n2-1}= v_{\frac n2-1}$ and $\hat v_{\frac n2+1}= v_{\frac n2+1}$.

Let us compute the area of this polygon defined by the vertices $\hat v_i$ which corresponds to the $n-$gon defined by vertices $v_i$ and where a substitution $\sigma$ is applied:
$$
\hat A_n=\frac{1}{2}\ds\sum_{i=1}^{n-2}(\hat y_i\hat x_{i+1}-\hat x_i\hat y_{i+1}).
$$ 
By remarking that the area for the quadrilateral defined by vertices $\hat v_n, \hat v_{\frac n2-1}, \hat v_{\frac n2}, \hat v_{\frac n2+1}$ is equal to $x_1$,
one has:
$$
\hat A_n=x_1+\frac{1}{2}\left(\ds\sum_{i=1}^{\frac n2-2}(\hat y_i\hat x_{i+1}-\hat x_i\hat y_{i+1})+\ds\sum_{i=\frac 
n2+1}^{n-2}(\hat y_i\hat x_{i+1}-\hat x_i\hat y_{i+1})\right).
$$
Because $\hat x_i=-\overline{x}_{n-i}$ and $\hat y_i=\overline{y}_{n-i}$ for all $i=1,\ldots, \frac n2 -2, \frac n2 +2, \ldots, n-1$, one obtains:
\begin{eqnarray*}
\hat A_n&=&x_1+\frac{1}{2}\left(\ds\sum_{i=1}^{\frac n2-2}(-\overline{y}_{n-i}\overline{x}_{n-(i+1)}+\overline{x}_{n-i}\overline{y}_{n-(i+1)})\right.\\
&&\left.\hspace{1,2cm}+\ds\sum_{i=\frac
n2+1}^{n-2}(-\overline{y}_{n-i}\overline{x}_{n-(i+1)}+\overline{x}_{n-i}\overline{y}_{n-(i+1)})\right).
\end{eqnarray*}
By fixing $j=n-i-1$, one directly obtains $\hat A_n=A_n$
and the result follows.
\end{DEM}

This lemma is not a proof that the largest $n-$gons with even $n\geq 4$ have a symmetry axis passing through vertices $(0,0)$ and $(0,1)$ in the quadratic formulation given in (\ref{quad-prog}). Numerical experiments reported below indicate however that this is always the case. Hence, we can just conjecture it.

\begin{CONJ}\label{conj}
The largest small $n-$gon with even $n\geq 4$ has a symmetry axis corresponding to the pending edge in its optimal diameter graph configuration.
\end{CONJ}   

Note that this conjecture is proved for $n=4$ by Reinhardt \cite{Rei1922} and for $n=6$ by Yuan \cite{Yua2004}.
For the latter case, the demonstration was not given in \cite{Gra1975}, however Graham used this result to find his little hexagon.
Moreover, Hansen's little octagon found in \cite{Aud2002} is a further evidence that the conjecture may be true.

\section{GloptiPoly}
\label{sec:GloptiPoly}

In 2000, Lasserre \cite{lasserre2000} proposed to reformulate a nonconvex polynomial optimization problem (POP)
\begin{equation}\label{pop}
\begin{array}{ll}
\min & g_0(x) \\
\mathrm{s.t.} & x \in X = \{x \in {\mathbb R}^n \: :\: g_i(x) \geq 0, \: i=1,\ldots,m\}
\end{array}
\end{equation}
where $g_0(x),g_1(x),\ldots,g_m(x)$ are multivariate polynomials,
as a linear infinite-dimensional moment problem,
in turn truncated into a primal-dual linear semidefinite programming (SDP) problem
\begin{equation}\label{sdp}
\begin{array}{c@{\hspace{2cm}}c}
\begin{array}{ll}
\min & c^T {\mathbf x} \\
\mathrm{s.t.} & A{\mathbf x} = b \\
& {\mathbf x} \in K
\end{array} &
\begin{array}{ll}
\max & b^T {\mathbf y} \\
\mathrm{s.t.} & {\mathbf z} = c-A^T{\mathbf y}\\
& {\mathbf z} \in K.
\end{array} 
\end{array}
\end{equation}
Using results on flat extensions of moment matrices and representations of polynomials positive on
semialgebraic sets, it was shown that under some relatively mild assumptions (implying in
particular that set $X$ is bounded), solving nonconvex POP
(\ref{pop}) amounts to solving a sufficiently large linear hence convex SDP problem (\ref{sdp}).
In this problem, vector $\mathbf y$ contains the moments of a probability measure supported on $X$
\begin{equation}
\label{moments}
{\mathbf y}_{\alpha} = \int_X x^{\alpha} d\mu(x)
\end{equation}
where we used the multi-index notation $x^{\alpha} = x^{\alpha_1}_1 x^{\alpha_2}_2 \cdots x^{\alpha_n}_n$.
Solving problem (\ref{sdp}) then amounts to optimizing over such probability measures.
If problem (\ref{pop}) has a finite number of global optimizers, then the optimal
probability measure is a linear combination of Dirac measures at these global optimizers.

In practice, a hierarchy of embedded SDP relaxations of increasing size are solved gradually,
for moments of increasing orders.
Convergence and hence global optimality can be guaranteed by examining a certain rank
pattern in the moment matrix, a simple task of numerical linear algebra. A user-friendly Matlab interface
called GloptiPoly was designed in 2002 to transform a given POP (\ref{pop})
into an SDP relaxation (\ref{sdp}) of given size in the hierarchy,
and then to call SeDuMi, a general-purpose conic solver \cite{gloptipoly2}. A new version 3 was released in 2007
to address generalized problem of moments, including POPs but also many other decision problems. The
interface was also extended to other public-domain conic solvers \cite{gloptipoly3}. Almost a decade after
the initial spark \cite{lasserre2000}, Lasserre summarized the theoretical and practical sides of the approach
in a monograph \cite{lasserre2009}. GloptiPoly is freely available for download at
\begin{verbatim}
homepages.laas.fr/henrion/software/gloptipoly
\end{verbatim}

For the case of the octagon ($n=8$) here is the GloptiPoly 3 Matlab code which is used to model the first SDP relaxation
in the hierarchy:

\begin{verbatim}
mpol x1 x2 x3 x4 x5 y1 y2 y3 y4 y5

f = 1/2*((x2+x3-4*x1)*y1+(3*x1-2*x3+x5)*y2+(3*x1-2*x2+x4)*y3+...
    (x3-2*x1)*y4+(x2-2*x1)*y5)+x1;

K = [(x1-x2)^2+(y1-y2)^2<=1, (-x1+x3-x5)^2+(y1-y3-y5)^2<=1, ...
(x1-x2+x4)^2+(y1-y2-y4)^2<=1, (-x1+x3)^2+(y1-y3)^2<=1, ...
(2*x1-x2-x3+x5)^2+(-y2+y3-y5)^2<=1, (2*x1-x2)^2+y2^2<=1, ...
(x1-x2)^2+(y1-y2-1)^2<=1, (2*x1-x2-x3)^2+(-y2+y3)^2<=1, ...
(x3-x5)^2+(-y3+y5)^2<=1, (-x1+x3-x5)^2+(y1-y3-y5+1)^2<=1, ...
(2*x1-x3+x5)^2+(-y3+y5)^2<=1, ...
(2*x1-x2-x3+x4+x5)^2+(-y2+y3+y4-y5)^2<=1,...
(-2*x1+x2-x4)^2+(y2-y4)^2<=1, (x1-x2+x4)^2+(y1-y2+y4-1)^2<=1, ...
(x1-x3)^2+(1-y1+y3)^2<=1, (x2-x4)^2+(y2-y4)^2<=1, ...
(2*x1-x3)^2+y3^2<=1, (2*x1-x2-x3+x4)^2+(-y2+y3+y4)^2<=1,...
x2-x3>=0, y1>=0, y2>=0, y3>=0, y4>=0, y5>=0, 0<=x1, x1<=1/2,
0<=x2, x2<=1, 0<=x3, x3<=1, 0<=x4, x4<=1, 0<=x5, x5<=1,...
y5^2==1-x5^2, y4^2==1-x4^2, y3^2==1-x3^2, y2^2==1-x2^2, ...
y1^2==1-x1^2];

P = msdp(max(f),K);
\end{verbatim}
These commands generate a GloptiPoly object {\tt P} which encodes the first
SDP relaxation (\ref{sdp}) of the nonconvex POP (\ref{pop}). Note that the 5 equality constraints
liking the $x_i$ and $y_i$ variables are exploited by the {\tt msdp} command
to reduce the size of the resulting SDP relaxation, using moment substitutions
as explained in \cite{gloptipoly3}. The SDP relaxation in primal-dual form (\ref{sdp})
can be retrieved with the following command:
\begin{verbatim}
[A,b,c,K] = msedumi(P);
\end{verbatim}
The SDP relaxation can be solved with the following command calling
the default conic solver (SeDuMi 1.3 in our case):
\begin{verbatim}
[stat,obj] = msol(P);
\end{verbatim}
We obtain {\tt stat = 0} indicating that GloptiPoly provides in output argument {\tt obj}
a valid upper bound on the global maximum, equal to $0.72686848$ (with 8 significant
digits). In this case we have to go deeper in the hierarchy and with the following
instructions we build and solve the second SDP relaxation:
\begin{verbatim}
P = msdp(max(f),K,2);
[stat,obj] = msol(P);
\end{verbatim}
We obtain {\tt stat = 1} indicating that GloptiPoly certifies numerically global optimality
(the moment matrix has approximately rank one), and it provides
in output argument {\tt obj} an upper bound $0.72686848$. With the command
\begin{verbatim}
double([x1 x2 x3 x4 x5])
\end{verbatim}
we can retrieve the solution (with 8 significant digits)
$x_1=0.26214172$, $x_2=0.67123417$, $x_3=0.67123381$, $x_4=0.90909242$, $x_5=0.90909213$.
This SDP problem is solved by SeDuMi in less than 5 seconds.
The quadratic objective function evaluated at the above solution is the same as
the computed upper bound to 11 significant digits. The
symmetry considerations of Lemma \ref{group} and Conjecture 
\ref{conj} indicate that $x_2=x_3$ and $x_4=x_5$ at the optimum,
and we see that the above solution achieves this to 5 digits for $x_2$
and to 6 digits for $x_4$.

Moreover, these results can be rigorously guaranteed by using Jansson's VSDP package which
uses SDP jointly with interval arithmetic \cite{vsdp}. The solution of an SDP problem can be
guaranteed at the price of solving a certain number of SDP problems of the
same size. In our case, VSDP solved 8 instances of the second SDP relaxation
to provide the guaranteed lower bound $0.72686845$ and guaranteed upper bound
$0.72686849$ on the objective function, namely the area of the octagon.

Generally speaking, note that the current status of SDP problem solving
is rather disappointing, in the sense that, to the best of our knowledge,
there is currently no backward stable SDP solver, and there is no efficient estimate
of conditioning of an SDP problem (by efficient we mean computationally less efficient
than solving the SDP problem itself). The VSDP package addresses indirectly these
issues in the sense that it uses interval arithmetic to provide guaranteed bounds
on the primal-dual solutions of problem (\ref{sdp}), but at the price of
a significant increase of the computational burden. In the lack of backward
stability guarantees for SDP solvers and efficient estimates of SDP problem conditioning,
we are not aware of a cheaper alternative to rigorously certify solutions
of SDP relaxations of nonconvex POPs.

\section{Numerical experiments}
\label{sec:NUM}

We applied GloptiPoly 3 and SeDuMi 1.1R3 to solve the nonconvex quadratic optimization problems.
In order to obtain accurate solutions, we let SeDuMi minimize the duality gap as much as possible.
We also tightened the tolerance parameters used by GloptiPoly to detect global optimality and extract
globally optimal solutions.  We used a 32 bit desktop personal computer with a standard configuration
and we report our numerical results to 8 significant digits.

\subsection{The largest small decagon}

\begin{figure}[ht]
\setlength{\unitlength}{6cm}
\begin{picture}(-3.3,1)(-0.9,0)%(40,22)(0,-20)
\thicklines
 \dottedline{0.001}(-0.3376, 0.1414)(0.2110, 0.9775)
 \dottedline{0.001}(-0.3376, 0.1414)(0.4453,  0.7635)
 \dottedline{0.001}(-0.5000, 0.4373)(0.4453,  0.7635)
 \dottedline{0.001}(-0.5000, 0.4373)(0.5, 0.4373)
 \dottedline{0.001}(-0.4453,  0.7635)(0.5000, 0.4373)
 \dottedline{0.001}(-0.4453,  0.7635)(0.3376, 0.1414)
 \dottedline{0.001}(-0.2110, 0.9775)(0.3376, 0.1414)
 \dottedline{0.001}(-0.2110, 0.9775)(0,0)
 \dottedline{0.001}(0,1)(0,0)
 \dottedline{0.001}(0.2110, 0.9775)(0,0)
  \put(-0.3376, 0.1414){\makebox(0,0){$\bullet$}}
 \put(-0.4453,  0.7635){\makebox(0,0){$\bullet$}}
 \put(-0.5000, 0.4373){\makebox(0,0){$\bullet$}}
 \put(-0.2110, 0.9775){\makebox(0,0){$\bullet$}}
 \put(0, 1){\makebox(0,0){$\bullet$}}
 \put(0.2110, 0.9775){\makebox(0,0){$\bullet$}}
 \put(0.4453, 0.7635){\makebox(0,0){$\bullet$}}
 \put(0.5000, 0.4373){\makebox(0,0){$\bullet$}}
 \put(0.3376, 0.1414){\makebox(0,0){$\bullet$}}
 \put(0,0){\makebox(0,0){$\bullet$}}
\dottedline{0.04}(-0.3376, 0.1414)(-0.5000, 0.4373)(-0.4453,  0.7635)(-0.2110, 0.9775)(0, 1)(0.2110, 0.9775)(0.4453, 0.7635)(0.5000, 0.4373)(0.3376, 0.1414)(0,0)(-0.3376, 0.1414)
\end{picture}
\caption{Largest Small Decagon.\label{Fig-decagon}}
\end{figure}
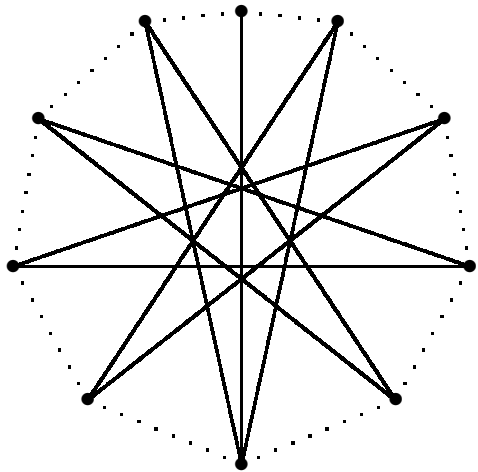

In the case $n=10$, we obtain the solution
$x_1=0.21101191$, $x_2=0.54864468$, $x_3=0.54864311$, $x_4=0.78292524$, $x_5=0.78292347$,
$x_6=0.94529290$, $x_7=0.94529183$ whose global optimality is guaranteed numerically
at the second SDP relaxation. This SDP problem, containing
2240 variables (size of the moment vector) and a semidefinite cone
of size 113 (size of the moment matrix),
is solved by SeDuMi in a little bit
more than 1 minute.
The objective function of the SDP relaxation, an upper bound on the exact global optimum,
is equal to $0.74913736$.  The quadratic
objective function evaluated at the above solution is the same to 10 significant digits.
The solution for the optimal decagon is drawn in Figure~\ref{Fig-decagon}.
Consistently with Conjecture 
\ref{conj}, we observe a symmetry axis on the optimal solution, namely $x_2=x_3$,
$x_4=x_5$ and $x_6=x_7$ up to 5 significant digits.

\subsection{The largest small dodecagon}

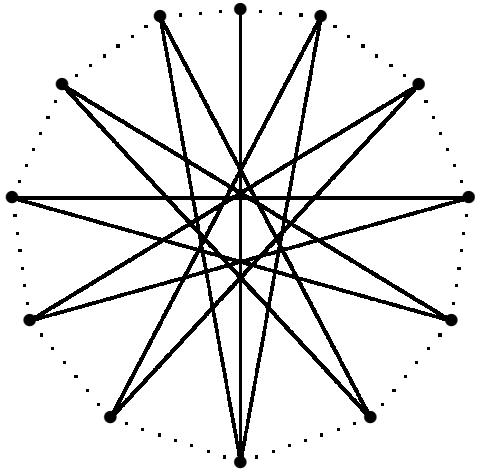
\begin{figure}[ht]
\setlength{\unitlength}{6cm}
\begin{picture}(-3.3,1)(-0.9,0)%(40,22)(0,-20)
\thicklines
 \dottedline{0.001}(-0.2853,    0.0972)(0.1762,    0.9844)
 \dottedline{0.001}(-0.2853,    0.0972)(0.3909,    0.8339)
 \dottedline{0.001}(-0.4623,    0.3123)(0.3909,    0.8339)
 \dottedline{0.001}(-0.4623,    0.3123)(0.5000,    0.5843)
 \dottedline{0.001}(-0.5000,    0.5843)(0.5000,    0.5843)
 \dottedline{0.001}(-0.5000,    0.5843)(0.4623,    0.3123)
 \dottedline{0.001}(-0.3909,    0.8339)(0.4623,    0.3123)
 \dottedline{0.001}(-0.3909,    0.8339)(0.2853,    0.0972)
 \dottedline{0.001}(-0.1762,    0.9844)(0.2853,    0.0972)
 \dottedline{0.001}(-0.1762,    0.9844)(0,0)
 \dottedline{0.001}(0,1)(0,0)
 \dottedline{0.001}(0.1762,    0.9844)(0,0)
 \put(-0.2853,    0.0972){\makebox(0,0){$\bullet$}}
 \put(-0.4623,    0.3123){\makebox(0,0){$\bullet$}}
 \put(-0.5000,    0.5843){\makebox(0,0){$\bullet$}}
 \put(-0.3909,    0.8339){\makebox(0,0){$\bullet$}}
 \put(-0.1762,    0.9844){\makebox(0,0){$\bullet$}}
 \put(0, 1){\makebox(0,0){$\bullet$}}
 \put(0.1762,    0.9844){\makebox(0,0){$\bullet$}}
 \put(0.2853,    0.0972){\makebox(0,0){$\bullet$}}
 \put(0.4623,    0.3123){\makebox(0,0){$\bullet$}}
 \put(0.5000,    0.5843){\makebox(0,0){$\bullet$}}
 \put(0.3909,    0.8339){\makebox(0,0){$\bullet$}}
 \put(0,0){\makebox(0,0){$\bullet$}}
\dottedline{0.04}(-0.2853, 0.0972)(-0.4623, 0.3123)(-0.5000, 0.5843)(-0.3909, 0.8339)(-0.1762, 0.9844)(0, 1.0000)(0.1762, 0.9844)(0.3909, 0.8339)(0.5000, 0.5843)(0.4623, 0.3123)(0.2853, 0.0972)(0, 0)(-0.2853, 0.0972)
\end{picture}
\caption{Largest Small Dodecagon.\label{Fig-dodecagon}}

\end{figure}

For $n=12$, without exploiting symmetry, the second SDP relaxation
contains 5640 variables (size of the moment vector) and a semidefinite cone of size 181
(size of the moment matrix). Such SDP problems are currently challenging for conic solvers,
even though recent progress on projection/regularization algorithms \cite{hm},
indicate that we may soon be able to solve routinely problems with semidefinite matrices
of size by the thousands. After approximately 25 minutes of CPU time, we obtain the following solution:
$x_1=0.17616131$, $x_2=0.46150224$, $x_3=0.46150519$, $x_4=0.67623091$, $x_5=0.67623301$, $x_6=0.85320300$,
$x_7=0.85320328$, $x_8=0.96231370$, $x_9=0.96231344$ featuring the expected symmetry
of Conjecture \ref{conj}.
The objective function is equal to $0.76072988$. The solution for the optimal dodecagon is drawn in Figure~\ref{Fig-dodecagon}.

\subsection{Exploiting the symmetry axis}

Remarking that Conjecture~\ref{conj} was formally proved for cases $n=4$ and $6$, \cite{Gra1975, Yua2004} and moreover that it is shown numerically for cases $n=8, 10$ and $12$ , we performed the determination of the following largest small polygons for $n=10$ to $16$ (with $n$ even) using this hypothesis of symmetry.

\begin{itemize}
\item Decagon ($n=10$): 

Recall from Lemma \ref{group} that problem (\ref{quad-prog}) is invariant
under the action of permutation group $G$. The SDP relaxations (\ref{sdp})
of the corresponding problem (\ref{pop}) are also invariant w.r.t. $G$.
It follows that several moments (\ref{moments}) are equal, for example
if $x_2$ and $x_3$ can be permuted we have
\[
\int_X x^{\alpha_2}_2 x^{\alpha_3}_3 d\mu(x) = \int_X x^{\alpha_3}_2 x^{\alpha_2}_3 d\mu(x) 
\]
for all integers $\alpha_2$ and $\alpha_3$. Therefore
vector $\mathbf y$ contains a lot of redundant entries that can
be removed. Similarly, the constraints in problem (\ref{sdp})
can be reduced significantly. See \cite{theobald} for a description
of how symmetry can be exploited in SDP relaxations of polynomial
optimization problems, and see also \cite{deklerk} for
a recent short survey on exploiting special structure in SDP problems,
in particular symmetry and invariance under a group of substitutions.

Rather than implementing the substitutions in the vector of moments
as described above, we went further and substituted all the variables
that can be permuted, e.g. we replaced $x_3$ with $x_2$, $x_5$ with $x_4$
and $x_7$ with $x_6$, respectively. The resulting SDP relaxations (\ref{sdp})
are thus significantly smaller. 
For comparison, the second SDP relaxation of the decagon quadratic problem
without substitutions
contains 2240 variables (size of the moment vector) and a largest semidefinite cone of size 113
(size of the moment matrix)
whereas the second SDP relaxation with substitutions contains 320 variables
and a largest semidefinite cone of size 41, a significant reduction.
This latter SDP problem is solved by SeDuMi in about 2 seconds, and
GloptiPoly can certify numerically global optimality of the solution
with a rank-one moment matrix. The objective function is equal to
$0.74913735$, and the solution is $x_1=0.21101121$, $x_2=0.54864181$,
$x_4=0.78292327$, $x_6=0.94529267$, consistently with the above
solution obtained without exploiting symmetry.

We also have tried to use VSDP jointly with SeDuMi 1.3 to certify rigorously the optimal decagon found
without exploiting symmetry, as we did for the octagon in Section 
\ref{sec:GloptiPoly}, but we could not obtain any meaningful result.
We have not tried to use alternative conic solvers such as e.g. SDPT3.
However, when solving the reduced second SDP relaxation, VSDP
provides the guaranteed lower bound $0.74913721$ and guaranteed upper
bound $0.74913740$ on the area of the optimal symmetric decagon.
The lower bound is also guaranteed for the nonsymmetric case.

\item Dodecagon ($n=12$): 

By exploiting symmetry, the second SDP
relaxation contains 680 variables and a semidefinite cone of size 61. It is solved after
about 8 seconds with SeDuMi, and the returned solution is as follows: $x_1=0.17616079$,
$x_2=x_3=0.46150096$, $x_4=x_5=0.67622897$, $x_6=x_7=0.85319926$, $x_8=x_9=0.96231045$
for an objective function equal to $0.76072986$.

Running VSDP jointly with SeDuMi 1.3 on this problem provides only a guaranteed
upper bound on the objective function, equal to $0.76072997$.

\item Tetradecagon ($n=14$):

GloptiPoly finds in $23$ seconds the largest small symmetric tetradecagon: $x_1=0.15100047, x_2=x_3=0.39733106,
x_4=x_5=0.59117050, x_6=x_7=0.76441599, x_8=x_9=0.89237421,x_{10}=x_{11}=0.97279813$ achieving
the objective function $0.76753100$. In \cite{Aud2007, Mos2005, Rei1922}, an upper bound on
the area of a small $n-$gon is given:
$$
A_n\leq \sin^2\left(\frac\pi{2n}\right)\cot\left(\frac\pi{n}\right).
$$ 
Hence, one has the following inequalities:
$$
0.76753100\leq A_{14}^*\leq 0.76893595
$$
on the area of the largest small tetradecagon.

\item Hexadecagon ($n=16$):

GloptiPoly finds in $276$ seconds the largest small symmetric hexadecagon: $x_1=0.13204787, x_2=x_3=0.34840959,
x_4=x_5=0.52343183, x_6=x_7=0.68719098, x_8=x_9=0.81912908, x_{10}=x_{11}=0.91836386,
x_{12}=x_{13}=0.97935563$  achieving
the objective function $0.77185969$.

Considering the above upper bound, one has
the following inequalities
$$
0.77185969\leq A_{16}^*\leq 0.77279135
$$
on the area of the largest small hexadecagon.

\end{itemize}

\section{Conclusion}
\label{CONC}

GloptiPoly can be efficiently used to find some largest small polygons with an even number $n$ of vertices. The octagon case ($n=8$)
is most efficiently solved than in \cite{Aud2002}: (i)~the optimal solution is now certified with 7 digits (using VSDP jointly with  interval arithmetic) and (ii)~the required CPU time is about 5 seconds instead of 100 hours (in 1997). Furthermore, the next open instance for the decagon ($n=10$) is solved using GloptiPoly in approximately
1 minute, and the dodecagon ($n=12)$ is solved in approximately 25 minutes. Note however that these solutions could not be
certified rigorously with VSDP and interval arithmetic.

Symmetry of the problem can be exploited to reduce further the dimension of the SDP relaxations and hence the accuracy of the results, even though
we cannot prove theoretically that solving the reduced problem is equivalent to solving the non-reduced problem, see Conjecture 
\ref{conj}. We just observe experimentally for small size instances that our conjecture is true for $n=8, 10, 12$. Note that
this conjecture was formally proved for cases $n=4, 6$ elsewhere in the technical literature. Moreover,
we provide the solutions for the largest small symmetric tetradecagon ($n=14$) and hexadecagon ($n=16$)
which are conjectured to be also the optimal nonsymmetric ones.
These numerical experiments tend to show that it seems to be also possible to solve the next open cases $n=18$ and $n=20$ 
as soon as symmetry is exploited. 

We note also that
these nonconvex quadratic problems are always solved globally at the second SDP relaxation, a phenomenon that we also observed
for many quadratic problems from the technical literature \cite{gloptipoly2}.

In future works, we have to certify and guarantee the solution obtained for the cases $n=10$ and $n=12$.
Finally, we are currently investigating further applications of this technique to
other nonconvex polynomial optimization problems arising in geometry.

%\end{article}
\end{document}